\documentclass[12pt]{article}

\usepackage{theorem}
\usepackage{amssymb}
\usepackage{amsmath}
\usepackage{amsfonts}
\usepackage{times}
\usepackage{graphicx}

\usepackage[paper=a4paper, top=2.5cm, bottom=2.5cm, left=2.5cm, right=2.5cm]{geometry}

\def \RR {\mathbb R}

\def \eps {\varepsilon}
\def \vphi {\varphi}

\newtheorem{theorem}{Theorem}[section]
\newtheorem{lemma}[theorem]{Lemma}

\newtheorem{proposition}[theorem]{Proposition}
\newtheorem{corollary}[theorem]{Corollary}
 {\theorembodyfont{\rmfamily}}
\newtheorem{definition}[theorem]{Definition}
\newtheorem{remark}[theorem]{Remark}

\def\myffrac#1#2 in #3{\raise 2.6pt\hbox{$#3 #1$}\mkern-1.5mu\raise 0.8pt\hbox{$
#3/$}\mkern-1.1mu\lower 1.5pt\hbox{$#3 #2$}}

\def\qed{\hfill $\vcenter{\hrule height .3mm
\hbox {\vrule width .3mm height 2.1mm \kern 2mm \vrule width .3mm
height 2.1mm} \hrule height .3mm}$ \bigskip}

\begin{document}

\title{Concentration of measures supported on the cube}
\author{Bo`az Klartag\thanks{School of Mathematical Sciences, Tel-Aviv
University, Tel Aviv 69978, Israel.
Email: klartagb@tau.ac.il}}
\date{}
\maketitle

\abstract{We prove a log-Sobolev inequality for a certain class of
log-concave measures in high dimension. These are the probability
measures supported on the unit cube $[0,1]^n \subset \RR^n$ whose
density takes the form $\exp(-\psi)$ where the function $\psi$ is
assumed to be convex (but not strictly convex) with bounded pure
second derivatives. Our argument relies on a transportation-cost
inequality \'a la Talagrand.  }

\section{Introduction}

Consider a cube $Q \subset \RR^n$ of sidelength $\ell$ parallel to
the axes, that is, $Q$ is a translation of the set $(0, \ell)^n
\subset \RR^n$ (or of its closure, equivalently). In this paper we
prove a concentration inequality for a class of probability measures
supported on $Q$.

\medskip  Write $| \cdot |$ for the standard Euclidean norm in $\RR^n$
and $B^n = \{ x \in \RR^n ; |x| \leq 1 \}$ is the
Euclidean unit ball centered at the origin. For a subset $A \subset
\RR^n$ denote $A + \eps B^n  = \{ x + \eps y; x \in A, y \in B^n
\}$, the $\eps$-neighborhood of the set $A$.

\begin{theorem} Let $\ell> 0, M \geq 0$ and let $Q \subset \RR^n$ be a cube of
sidelength $\ell$ parallel to the axes. Let $\mu$ be a probability measure supported on $Q$
with density $\exp(-\psi)$ for a convex function $\psi: Q \rightarrow
\RR$ such that
\begin{equation} \partial^{ii} \psi(x) \leq M \quad \quad \quad \text{for all} \ x \in Q, i=1,\ldots,n.
\label{eq_2123} \end{equation}  Suppose that $A \subseteq \RR^n$ is a measurable set
with $\mu(A) \geq 1/2$. Then, for all $t > 0$,
\begin{equation}
 \mu \left( A + t B^n \right) \geq 1 - \exp \left(-
t^2 / \alpha^2 \right)   \label{eq_2126}
\end{equation}
where $\alpha = \alpha(\ell,M) = 3 \ell e^{M \ell^2/8}$.
\label{thm_2140}
\end{theorem}

Theorem \ref{thm_2140} is equivalent to a logarithmic Sobolev
inequality and to a concentration inequality for Lipschitz
functions, see Section \ref{sec4} below. In probabilistic
terminology, we consider  uniformly bounded random variables
$X_1,\ldots,X_n$, possibly dependent, whose joint distribution
satisfies the convexity/concavity assumptions of Theorem
\ref{thm_2140}. Our results correspond to bounds for the variance
and tail distribution
 of $f(X_1,\ldots,X_n)$ where $f$ is a Lipschitz function on $\RR^n$.

\medskip We emphasize that we are not assuming any product structure, any symmetries
nor strict convexity for the function $\psi$ from Theorem
\ref{thm_2140}. There is a vast body of literature pertaining to the
case in which the measure $\mu$ is an arbitrary product measure in
the cube, see Talagrand \cite{talagrand_ihes}, Marton \cite{marton},
Dembo and Zeitouni \cite{DZ}, Ledoux \cite{ledoux_product} and
others. When the function $\psi$ from Theorem \ref{thm_2140} admits
a uniform positive lower bound for the Hessian, the conclusion of
Theorem \ref{thm_2140} is well-known and essentially goes back to
Bakry and \'Emery \cite{BE}.

\medskip How
can we produce  probability measures satisfying the assumptions of Theorem
\ref{thm_2140} with, say, $M = 1$? One may begin with the standard
Gaussian density in $\RR^n$, the function $$ \gamma_n(x) = (2
\pi)^{-(n/2)} \exp(-|x|^2/2) \quad \quad \quad (x \in \RR^n). $$ The
restriction of $\gamma_n$ to any cube $Q \subset \RR^n$, normalized to be a probability density, surely
satisfies the assumptions of Theorem \ref{thm_2140} with $M = 1$.
Furthermore, begin with any  probability density $\rho:
\RR^n \rightarrow [0, \infty)$ which is log-concave (that is, the function $-\log \rho$ is convex). Consider the convolution
$$ f(x) = (\rho * \gamma_n)(x) = \int_{\RR^n} \rho(y) \gamma_n(x-y) dy. $$
Then $f$ is a smooth, log-concave probability density according to the Pr\'ekopa-Leindler inequality.
Furthermore, it is straightforward to verify that for any $x \in \RR^n$,
\begin{equation}  (\nabla^2 \log f)(x) \geq -Id
\label{eq_2146}  \end{equation}
in the sense of symmetric matrices, where $Id$ is the identity matrix and $\nabla^2 \log f$ is the Hessian of $\log f$.
We conclude that the probability measure on the cube $Q$ whose density is proportional to the restriction of $f$ to $Q$,
satisfies the assumptions of Theorem \ref{thm_2140} with $M = 1$.
It is also possible to view
the probability densities that appear in Theorem \ref{thm_2140}
as convex perturbations of probability densities
proportional to $x \mapsto \exp(x \cdot v)$ on the cube. Here $x \cdot v$ is the standard scalar
product of $x,v \in \RR^n$.

\medskip
One cannot replace $\alpha(\ell,M)$ in Theorem
\ref{thm_2140} by a dimension-free expression that is
subexponential in $M \ell^2$, see Remark \ref{rem_1703} below.
We say that a vector $x \in \RR^n$ is proportional to one of the standard unit vectors
when it has at most one non-zero entry. A unit cube has sidelength one.
Theorem
\ref{thm_2140} will be deduced from the following  result:

\begin{theorem}
Let $R \geq 1$ and let $Q \subset \RR^n$ be a unit cube
parallel to the axes. Let $\mu$ be a probability measure supported on $Q$ with
a log-concave density $f$ such that
\begin{equation} f \left( \lambda x + (1 - \lambda) y \right) \leq R \left[  \lambda f(x) + (1 - \lambda) f(y) \right]
\label{eq_1635} \end{equation} for any $0 < \lambda < 1$ and any $x,
y \in Q$ for which $x - y$ is proportional to one of the standard
unit vectors. Suppose that $A \subseteq \RR^n$ is a measurable set with
$\mu(A) \geq 1/2$. Then for all $t > 0$,
$$
 \mu \left( A + t B^n \right) \geq 1 - \exp \left(-
t^2 / \alpha^2 \right) $$where $\alpha =
\alpha(R) = 3 R$. \label{cor_1636}
\end{theorem}

The inequality (\ref{eq_1635}) holds true with $R=1$ when $f$ is a
convex function. By degenerating Theorem \ref{cor_1636} to the petty
case where $R = 1$ we arrive at the following peculiar corollary:

\begin{corollary}
Let $Q \subset \RR^n$ be a unit cube.
Let $\mu$ be a probability measure on $Q$ whose density is both log-concave and convex in $Q$.
Then for any measurable $A \subseteq \RR^n$ and $t > 0$,
$$
 \mu(A) \geq 1/2 \quad \quad \Longrightarrow \quad \quad \mu \left( A + t B^n \right) \geq 1 - \exp \left(-
t^2 / 9 \right). $$ \label{cor_1543}
\end{corollary}

A moment of reflection reveals that there do exist positive, integrable functions on the cube that are simultaneously
log-concave and convex, such as $x \mapsto \left[ b + (x \cdot v) \right]^p$ for $p \geq 1$. It is also evident
that one cannot eliminate neither the log-concavity assumption nor the convexity assumption
from Corollary \ref{cor_1543}.

\medskip The proof of  Theorem \ref{cor_1636} uses transportation of measure  in
order to analyze the deficit in the Pr\'ekopa-Leindler inequality,
an idea proposed also in Eldan and Klartag \cite{EK}. Rather than
working directly with the supremum-convolution, we prefer to analyze
another expression that somewhat resembles the relative-entropy
functional. Let us  shed
some light on this expression. Suppose that $f$ and $g$ are
non-negative functions defined on $\RR^n$.  For a point $x \in
\RR^n$ in which $f$ is positive and differentiable, and for a point
$y \in \RR^n$ in which $g$ is positive, we set
\begin{equation}
S_{y} \{g, f \}(x) = f(x) \log \frac{g(y)}{f(x)} - \nabla f(x) \cdot
(y - x).
 \label{eq_501}
\end{equation}
Denote $Supp(f) = \{ x ; f(x) \neq 0 \}$. For functions $f, g: \RR^n \rightarrow [0, \infty)$ and
 a map $T: Supp(f) \rightarrow Supp(g)$ abbreviate
 \begin{equation}
S_T \{ g, f \} (x) = S_{T(x)} \{ g, f \}(x) \quad \quad \quad \quad
(x \in Supp(f)), \label{eq_402}
\end{equation}
assuming that the right-hand side is well-defined.
Next, suppose that $f$ and $g$ have a finite, positive integral. A measurable map $T: \overline{Supp(f)} \rightarrow \overline{Supp(g)}$
is called a {\it transportation map from $f$ to $g$} if for any measurable set
$A \subseteq Supp(g)$,
$$ \left( \frac{1}{\int g} \right) \int_{A} g  = \left( \frac{1}{\int f} \right) \int_{T^{-1}(A)} f.  $$
That is, $T$ pushes forward the probability measure whose density is
proportional to $f$, to the probability measure whose density is
proportional to $g$. Two important examples of
transportation maps in $\RR^n$ are the Brenier map
\cite{brenier} and the Knothe map \cite{knothe}.

\begin{definition} Let $f,g$ be two non-negative functions on $\RR^n$ with a finite, positive integral. Assume that $f$
is differentiable almost-everywhere in $Supp(f)$. Set
\begin{equation}  Tire(g \, || \, f) = \sup_T \left[ \int_{Supp(f)} S_T \{g, f \}(x) dx \, - \, \left( \int f \right) \log \frac{ \int g }{\int f  }
\right], \label{eq_1046} \end{equation}
where the supremum  runs over all transportation maps $T$ from $f$ to $g$
for which the integral of $S_T \{ g, f \}$ is well-defined. Here,   Tire is an acronym  of ``Translation-Invariant
Relative Entropy''.
 \label{tire}
\end{definition}

The notion is indeed translation-invariant: For functions $f,g$ as in Definition \ref{tire} and for $x_0 \in \RR^n$,
denoting $\tau_{x_0}(g)(x) = g(x-x_0)$ we have
$$  Tire(g \, || \, f) \, - \, Tire( \tau_{x_0}(g) \, || \, f )  =
\int_{\RR^n} \left( \nabla f(x)
\cdot x_0 \right) dx = 0,  $$
where we assume that $f$ is locally-Lipschitz and vanishes at infinity in order to justify the integration by parts.
In the log-concave case, the quantity $Tire(g \, || \, f)$ is indeed related to relative entropy as is demonstrated in Lemma \ref{lem_1035} below.

\medskip The remainder of this paper is devoted to the proofs of the aforementioned theorems
and  to related results.
 We write $\overline{A}$ for the closure
of the set $A$, and $\log$ stands for the natural
logarithm. By ``measurable'' we always mean Borel-measurable.
Needless to say, it is certainly possible to consider $Tire(g \, || \, f )$ for non-negative functions
defined only on a subset of $\RR^n$ by treating such functions as zero outside their original domain of definition.

\medskip
{\it{Acknowledgements.}} I thank Dario Cordero-Erausquin and Ronen Eldan for interesting,
related discussions. I am grateful to Nathael Gozlan and to
the anonymous referee for their valuable suggestions
and for correcting a mistake
in an earlier version of this manuscript.
The research was supported in part
 by the Israel
Science Foundation and by a Marie Curie Reintegration Grant from the
Commission of the European Communities.

\section{Convex functions on an interval}

 Let $I,J
\subset \RR$ be two intervals of finite, positive length and let $f,g$ be positive, integrable
functions defined on $I,J$ respectively. The {\it monotone transportation map from $f$
to $g$} is the map $T: \overline{I} \rightarrow \overline{J}$ defined via
$$ \left( \frac{1}{\int_I f} \right) \int_I f(t) 1_{\{ t < x \}} dt = \left( \frac{1}{\int_J g} \right) \int_J
 g(t) 1_{ \{t < T(x) \}} dt \quad \quad \quad (x \in \overline{I}) $$ where $1_{\{ t < x \}}$ equals
one when $t < x$ and vanishes otherwise. The map $T$ is uniquely
defined, as $f,g$ are positive and integrable.
Furthermore, $T$ is
an absolutely-continuous, strictly-increasing function.
Observe that the monotone transportation in one
dimension is indeed a transportation map and that for almost every $x \in I$, \begin{equation}
 T^{\prime}(x) = \left( \frac{\int_J g}{\int_I f} \right) \frac{f(x)}{g(T(x))}.
\label{eq_1047} \end{equation}
We will frequently encounter the case where $I = J$.
Clearly, in this case  $T(x) = x$ for $x
\in \partial I$, where $\partial I$ are the two endpoints of the
interval $I$.
  Our goal in this section
is to prove the following transportation-cost inequality in one dimension:

\begin{proposition} Let $R \geq 1$ and let $I \subset \RR$ be an interval of length
one. Let $f: I \rightarrow (0, \infty)$ be an absolutely-continuous function which
satisfies
$$ f(\lambda x + (1-\lambda) y) \leq R \left[ \lambda f(x) + (1-\lambda) f(y) \right] \quad \quad \text{for all} \
\ x,y \in I, 0 < \lambda < 1. $$ Let $g$
be a positive, integrable function on $I$, and let $T$ be the
monotone transportation map from $f$ to $g$. Then,
\begin{align} \label{eq_147} \int_I \left| T(x) -x    \right
|^2 f(x) dx  & \leq C R^{2} \left[ \int_{I} S_T \{ g, f \} -
\left(
\int_{I} f \right) \log \frac{\int_{I} g}{\int_{I} f} \right] \\
& \leq C R^{2} \cdot Tire( g \, || \, f ), \nonumber
\end{align}
where $C \leq 40/9$ is a universal constant.
\label{prop_142}
\end{proposition}

The proof of Proposition \ref{prop_142} requires a few
lemmata. Our first lemma is essentially an infinitesimal version of the
Pr\'ekopa-Leindler inequality, and its proof follows the
transportation proofs given by Barthe \cite{barthe},
Cordero-Erausquin \cite{DCE}, Henstock-Macbeath \cite{HM} and Talagrand \cite{talagrand}.
For $t \in \RR$ denote
$$ \Lambda(t) = \min \{ |t|, t^2 \}. $$

\begin{lemma} Let $I \subset \RR$ be an interval of  finite, positive length. Let $f,g$ be
positive, integrable functions on $I$ with $f$ being absolutely
continuous. Then,
$$ \int_I \Lambda
\left( T^{\prime}(x) - 1 \right) f(x) dx
\, \leq \, \frac{10}{3} \left[ \int_{I} S_T \{ g, f \} - \left( \int_{I} f \right) \log
\frac{\int_{I} g}{\int_{I} f} \right]
$$ where $T$ is the monotone transportation map from $f$ to $g$.
\label{lem_1145}
\end{lemma}

\textbf{Proof:} We use (\ref{eq_1047}) and compute
\begin{align} \nonumber
\int_{I} & S_T \{ g, f \}  = \int_{I} \left[ f(x) \log \frac{g(
T(x))}{f(x)} -  f^{\prime}(x)  (T(x) - x) \right] dx \\
\nonumber   & = \left( \int f \right) \log \frac{\int g}{\int f} +
\int_{I} \left[ f(x) \log \frac{1}{T^{\prime}(x)} - f^{\prime}(x)
(T(x) - x) \right] dx \\ &= \left( \int f \right) \log \frac{\int
g}{\int f} + \int_{I} \left[ -f(x) \log {T^{\prime} (x)} +
f(x) (T^{\prime}(x) - 1) \right] dx  \nonumber
\end{align}
where the integration by parts is legitimate as $f(x) (T(x) - x)$ is
an absolutely-continuous function in $I$ that vanishes on
$\partial I$. In order to complete the proof of the lemma it remains
 to show that for all $x > 0$,
\begin{equation}
- \log x + (x - 1)  \geq  \frac{3}{10} \cdot \min \{ |x-1|, (x-1)^2 \}. \label{eq_1124}
 \end{equation}
Indeed, for $0 < x \leq 2$ we use the Cauchy-Schwartz inequality and obtain
\begin{align*}
 (x -1) & - \log x = \int_1^x \frac{t-1}{t} dt = \int_0^{|x-1|}
\frac{t}{1 + sgn(x-1) t} dt \geq \int_0^{|x-1|}
\frac{t}{1 + t} dt \\ & \geq \left. \left( \int_0^{|x-1|} t dt \right)^2 \right / \left( \int_0^{|x-1|} (1 + t) t dt \right)
= \frac{(x-1)^2}{2(1 + 2 |x-1| /3)} \geq \frac{3 (x-1)^2}{10},
\end{align*}
where $sgn(x) = 1$ for $x > 0$ and $sgn(x) = -1$ for $x < 0$.
The inequality (\ref{eq_1124}) is valid in particular for $x
= 2$. For $x > 2$ the derivative of the left-hand side in (\ref{eq_1124}) exceeds
that of the right-hand side. Hence (\ref{eq_1124}) holds true for all $x
> 0$.
\qed

\begin{remark} {\rm The proof of Lemma \ref{lem_1145} admits a generalization
 to $n$ dimensions, in which one
utilizes the Brenier map in place of the transportation map $T$. See
Barthe \cite{barthe2} and McCann \cite{mccann} for related
arguments. In this way one obtains the inequality
\begin{equation}
Tire(g \, || \,
f) \geq 0, \label{eq_1831}
\end{equation}
which is valid for any Lipschitz, non-negative, compactly-supported
functions $f$ and $g$ on $\RR^n$ with a finite, positive integral.
\label{rem_1809}}
\end{remark}

\begin{lemma} Let $R \geq 1$ and let $I \subset \RR$ be an interval of length one.
Assume that $\rho$ is a positive, integrable function on $I$ that satisfies
\begin{equation} \rho(\lambda x + (1-\lambda) y) \leq R \left[ \lambda \rho(x) + (1-\lambda) \rho(y) \right] \quad \text{for all} \
\ x,y \in I, 0 < \lambda < 1. \label{eq_2106} \end{equation} Then for any $a,b \in I$ with
$a < b$,
$$  \int_{a}^b \rho(x) dx \leq  \frac{R}{2}  \left[ \rho(a) + \rho(b) \right].
$$
 \label{lem_244}
\end{lemma}

\textbf{Proof:} We simply integrate (\ref{eq_2106}) over $\lambda \in [0,1]$. Since $b-a \leq 1$, then
$$
 \int_a^b \rho \leq \int_0^1 \rho(\lambda a + (1 -
\lambda) b) d \lambda \leq R \int_0^1 \left[ \lambda \rho(a) +
(1-\lambda) \rho(b) \right] d \lambda = R \frac{\rho(a) +
\rho(b)}{2},
$$ and the lemma is proven.
\qed

\medskip
The following lemma is a one-dimensional Poincar\'e-type inequality.
The proof closely follows  the argument by Cheeger \cite{cheeger}.
Recall that $\Lambda(t) = \min \{ |t|, t^2 \}$.

\begin{lemma} Let $I \subset \RR$ be an interval of length
one and let $R \geq 1$. Let $\rho$ be a positive, integrable  function on $I$
with
$$ \rho(\lambda x + (1-\lambda) y) \leq R \left[ \lambda \rho(x) + (1-\lambda) \rho(y) \right] \quad \quad \text{for all} \
\ x,y \in I, 0 < \lambda < 1. $$
Then for any absolutely-continuous
function $f: \overline{I} \rightarrow \RR$ with $f|_{\partial I} =
0$,
\begin{equation} \int_I  \Lambda \left( f \right)  \rho \leq \frac{4}{3} R^{2}  \int_I \Lambda \left( f^{\prime}
\right)  \rho. \label{eq_846}
\end{equation} \label{lem_poincare1D} Here, $\partial I$ consists of
the two endpoints of the interval $I$.
 \label{lem_233}
\end{lemma}

\textbf{Proof:}  Multiplying $\rho$ by a constant, we may assume that $\int_I \rho = 1$. Let $g$ be an absolutely-continuous, non-negative
function with $g|_{\partial I} = 0$.
In the first
step of the proof we show that
\begin{equation} \int_I g \rho \leq \frac{R}{2}  \int_I \left| g^{\prime} \right| \rho.
\label{eq_847} \end{equation}   Denote $J = g(I) = \{ g(x) ; x
\in I \}$, an interval whose left boundary point is zero. We
apply the change of variables
 $y = g(x)$   and conclude that
 \begin{equation}  \int_I \left| g^{\prime} (x) \right| \rho(x) dx = \int_J
\left( \sum_{x \in g^{-1}(y)} \rho(x) \right) dy.
\label{eq_320} \end{equation}
See, e.g., Leoni \cite[Theorem 3.65]{leoni} for a proof of this change of variables formula.
For any $0 \neq y \in J$ consider the open
set $I_y = \{ x \in I ; g(x) > y \}$. When $y$ is a regular non-zero
value of $g$, the open set $I_y$ is a finite union of intervals with
disjoint closures. According to Lemma \ref{lem_244}, for any such
$y$,
$$ \int_{I_y} \rho \leq \frac{ R}{2} \left[ \sum_{x \in g^{-1}(y)} \rho(x) \right].
$$
The one-dimensional Sard's lemma for absolutely-continuous functions
(see, e.g., Leoni \cite[Remark 8.9]{leoni})
implies that almost any $y \in J$ is a regular value of $g$.
Therefore, from (\ref{eq_320}) we obtain
\begin{align*}  \int_I \left| g^{\prime} (x) \right| \rho(x) dx \geq  \frac{2}{R} \int_J
\left( \int_{\{ x ; g(x) > y \} } \rho(x) dx  \right) dy =
\frac{2}{R} \int_I g \rho
\end{align*}
where the last equality follows from  application of Fubini's
theorem. Thus (\ref{eq_847}) is proven. In order to prove
(\ref{eq_846}), observe that for any $x \geq 0$  and $0 \leq y \leq
1$,
\begin{equation} x y \leq R \Lambda(x) + \frac{y^2}{4 R}. \label{eq_209}
\end{equation}
Indeed, (\ref{eq_209}) holds for $x \geq 1$ since the coefficient in
front of $\Lambda(x)$ is at least one, and  (\ref{eq_209}) may be
directly proven for $0 \leq x \leq 1$  by completing a square. Let
$f: \overline{I} \rightarrow \RR$ be an absolutely-continuous
function with $f|_{\partial I} = 0$. Applying (\ref{eq_847}) with $g
= \Lambda(|f|)$ and using (\ref{eq_209}) we see that
\begin{align*} \int_I \Lambda(|f|) \rho  & \leq \frac{R}{2} \int_I  \Lambda^{\prime}(|f|) |f^{\prime}|
\rho \leq R \int_I \left| f^{\prime} \right| \min \{ |f|, 1
\} \rho \\ & \leq R \cdot \left[ R \int_I
\Lambda(|f^{\prime}|) \rho + \frac{1}{4 R} \int_I \min
\{ |f|^2, 1 \} \rho \right]
\\ & \leq R^2  \int_I
\Lambda(|f^{\prime}|) \rho + \frac{1}{4} \int_I \Lambda(|f|) \rho.
\end{align*}
By subtracting the right-most summand from  the  left-hand side, we
deduce (\ref{eq_846}).\qed

\textbf{Proof of Proposition \ref{prop_142}:} Since $T(x) = x$ for
$x \in \partial I$ we may invoke Lemma \ref{lem_233} and conclude
that
\begin{align*}
 \int_I \Lambda \left( T(x) - x \right) f(x) dx  & \leq \frac{4}{3} R^2
\int_I \Lambda \left(  T^{\prime}(x)
- 1 \right ) f(x) dx \\
& \leq \frac{40}{9} R^2  \left[ \int_{I} S_T \{ g, f \} - \left( \int_{I}
f \right) \log \frac{\int_{I} g}{\int_{I} f} \right]
\end{align*}
where  we used Lemma \ref{lem_1145} in the last passage. Since $I$ is an interval of length one
and $T: \overline{I} \rightarrow \overline{I}$, then for any $x \in I$ we have $|T(x) - x| \leq 1$. Consequently, for
any $x \in I$,
$$ \Lambda \left( T(x) - x \right) = \min \left \{ |T(x) - x|^2, |T(x) - x| \right \} = |T(x) - x|^2. $$
This completes the proof of (\ref{eq_147}). The proposition now follows from the definition of $Tire(g \, || \, f)$. \qed

\section{Induction on the dimension}

In this section we obtain higher-dimensional analogs of Proposition
\ref{prop_142}.

\begin{theorem}
Let $R \geq 1$ and let $Q \subset \RR^n$ be a unit cube  parallel to the axes.
Assume that $f: Q \rightarrow (0, \infty)$ is a Lipschitz function with
\begin{equation} f \left( \lambda x + (1 - \lambda) y \right) \leq R \left[ \lambda f(x) + (1-\lambda) f(y) \right] \label{eq_107} \end{equation} for any $0 < \lambda < 1$ and any $x,
y \in Q$ for which $x - y$ is proportional to one of the standard
unit vectors in $\RR^n$.
  Let $g$ be a positive,
integrable function on $Q$. Then there exists a
transportation map $T$ from $f$ to $g$ such that
\begin{align} \label{eq_1655} \int_Q \left| T(x) -x    \right
|^2 f(x) dx  & \leq C R^{2} \left[ \int_{Q} S_T \{ g, f \} -
\left(
\int_{Q} f \right) \log \frac{\int_{Q} g}{\int_{Q} f} \right] \\
& \leq C R^{2} \cdot Tire( g \, || \, f ), \nonumber
\end{align}
\label{thm_103}
 where $C \leq 40/9$ is a universal constant.
\end{theorem}

The requirement that $f$ be a Lipschitz function should not be taken
too seriously, as it may easily be replaced by other types of
regularity assumptions. Theorem \ref{thm_103} will be proven by
induction on the dimension, where the induction step is going to be
Proposition \ref{prop_142} in disguise. Throughout this section we
use
$$  x = (y,r) \in \RR^{n-1} \times \RR $$
as coordinates in $\RR^n$.  For a function $f$ defined on a subset
of $\RR^n$ and for $y \in \RR^{n-1}$, we write
$$ f_y(r) = f(y,r)  $$
whenever $(y,r)$ is in the domain of definition of $f$. Abbreviate
$\pi(y,r) = y$. For a subset  $A \subseteq \RR^n$  denote $\pi(A)
= \{ \pi(x) ; x \in A \}$. For a non-negative, integrable function
$f$ defined on a subset $A \subseteq \RR^n$,
 we set
$$ \pi(f)(y) = \int_{-\infty}^{\infty} f_y(r) 1_{ \{ (y,r) \in A \}} dr
\quad \quad \quad (y \in \pi(A)). $$
Let $K \subseteq \RR^n$ be a convex set. Let $f, g$ be positive,
integrable functions on $K$. We say that a map $T: K \rightarrow K$
 {\it transports the last coordinate monotonically} if there
exists a map $P: \pi(K) \rightarrow \pi(K)$ such that for almost any
$y \in \pi(K)$, the function $g_{P(y)}$ is integrable and
furthermore
\begin{equation}
 T(y, r) = (P(y), T_{y}(r))  \label{eq_1226} \end{equation}
for any $r$ with $(y,r) \in K$,
 where $T_{y}$ is the monotone transportation map from $f_y$ to $g_{P(y)}$.
The following lemma is a corollary to
Proposition \ref{prop_142}.

\begin{lemma} Let $R \geq 1$. Let $Q = A \times I \subset \RR^n$ where $I \subset \RR$ is an interval of length
one and $A \subset \RR^{n-1}$ is a convex set. Assume that $f$ is a positive,
Lipschitz  function on $Q$, and that
\begin{equation} f \left( \lambda x_1 + (1 - \lambda) x_2 \right) \leq R \left[ \lambda f(x_1) + (1 - \lambda) f(x_2) \right]
\label{eq_1527} \end{equation}
 for any $0 < \lambda < 1$ and any $x_1,
x_2 \in Q$ for which $x_1 - x_2$ is proportional to one of the standard
unit vectors. Let $g$ be a positive, integrable function on $Q$.
Assume that $T: \overline{Q} \rightarrow \overline{Q}$ is a measurable map that transports
the last coordinate monotonically. Then,
\begin{align} \label{eq_1239} \int_Q    \left |  T_y(r) -r \right |^2 f(y,r) dy dr   \, \leq \, C R^{2 }  \left[ \int_Q S_{T} \{ g,
f \} - \int_{\pi(Q)} S_{P} \{ \pi(g), \pi(f) \} \right]
\end{align}
 where $P$ and $T_y$ are as in (\ref{eq_1226}), and $C \leq 40/9$ is a universal constant.
 \label{lem_1232}
\end{lemma}

\textbf{Proof:} According to the definitions
(\ref{eq_501}) and (\ref{eq_402}), for almost any $(y,r) \in Q$,
 \begin{equation} S_{T} \{g, f \}(y,r)  = S_{T_y} \{
g_{P(y)}, f_y \}(r) - \nabla_y f(y, r) \cdot (P(y) - y) \label{eq_1124_}
\end{equation}
where $\nabla_y$ is the gradient in the $y$-variables.
Thanks to our assumptions on $f$ we may safely
differentiate under the integral sign, thus
\begin{equation} \nabla \pi(f)(y) \cdot (P(y) - y) = \int_I  \nabla_y f(y,r)  \cdot (P(y) -
y) dr \label{eq_1124__} \end{equation}
for almost any choice of $y$. From (\ref{eq_1124_}) and (\ref{eq_1124__}),
\begin{equation}
\int_{I} S_{T_y} \{ g_{P(y)}, f_y \}(r) dr = \int_I S_{T} \{g, f
\}(y,r) dr + \nabla \pi(f)(y) \cdot (P(y) - y) \label{eq_514}
\end{equation}
for almost any choice of $y$. The equality (\ref{eq_514}) may be reformulated
as
\begin{equation}
\int_{I} S_{T_y} \{ g_{P(y)}, f_y \} - \left( \int_{I} f_y \right) \log
\frac{\int_{I} g_{P(y)}}{\int_{I} f_y} = \int_I
S_{T} \{ g, f \} (y,r) dr - S_{P} \{ \pi(g), \pi(f) \}(y). \label{eq_448}
\end{equation}
We may apply Proposition \ref{prop_142} thanks to
our assumption (\ref{eq_1527}) and obtain that
\begin{equation}
\int_I  \left |  T_y(r) -r \right |^2 f_y(r) dr  \, \leq \,
C R^{2} \left[ \int_{I} S_{T_y} \{ g_{P(y)}, f_y \} - \left( \int_{I} f_y
\right) \log \frac{\int_{I} g_{P(y)}}{\int_{I} f_y} \right]. \label{eq_1442}
\end{equation}
By combining (\ref{eq_448}) and (\ref{eq_1442}) we see that for almost any $y \in
\pi(Q)$,
\begin{align} \label{eq_1240} \int_I  \left |  T_y(r) -r \right|^2 f_y(r) dr  \, \leq \, C R^{2}  \left[ \int_I S_{T} \{ g,
f \}(y,r) dr - S_{P} \{\pi(g), \pi(f) \}(y) \right].
\end{align}
We now integrate (\ref{eq_1240}) over $y \in \pi(Q)$ and deduce
(\ref{eq_1239}). \qed

\medskip \textbf{Remark:} The identity (\ref{eq_1124__}) is the only place
in the proof of Theorem \ref{thm_103} where we use the assumption
that $Q$ is a cube or  a box, rather than, say, a parallelepiped.

\begin{lemma} Let $R \geq 1$ and let $Q \subset \RR^n$ be a cube parallel to the axes. Assume that $f: Q \rightarrow \RR$ is
a Lipschitz  function on $Q$, such that
\begin{equation} f \left( \lambda x_1 + (1 - \lambda) x_2 \right) \leq R \left[ \lambda f(x_1) + (1-\lambda) f(x_2) \right]
\label{eq_146} \end{equation} for any $0 < \lambda < 1$ and any $x_1,
x_2 \in Q$ for which $x_1 - x_2$ is proportional to one of the standard
unit vectors in $\RR^n$.
Then also
$$ \pi(f) \left( \lambda y_1 + (1 - \lambda) y_2 \right) \leq R \left[ \lambda \pi(f)(y_1) + (1-\lambda) \pi(f)(y_2) \right] $$
for any $0 < \lambda < 1$ and any $y_1,
y_2 \in \pi(Q)$ for which $y_1 - y_2$ is proportional to one of the standard
unit vectors in $\RR^{n-1}$.
 \label{lem_2215}
\end{lemma}

\textbf{Proof:} Fix $i=1,\ldots,n-1$ and let $e_i$ be the $i^{th}$
standard unit vector. Condition (\ref{eq_146}) implies that for any
$y \in \RR^{n-1}, t_1, t_2, r \in \RR, 0 < \lambda < 1$ such that
$(y + t_1 e_i, r) \in Q$ and $(y + t_2 e_i, r) \in Q$,
$$ f \left( y + \left( \lambda t_1 + (1 - \lambda) t_2 \right)  e_i, r  \right) \leq R \left[ \lambda f \left( y + t_1 e_i, r \right) +
 (1-\lambda) f \left( y + t_2
e_i, r \right) \right].  $$
 Let $I$ be the interval for
which $Q = \pi(Q) \times I$. Integrating with respect to $r$ we have
\begin{align} \nonumber
 \pi( f) & \left( y + \left( \lambda t_1 + (1 - \lambda) t_2 \right) e_i \right)   = \, \int_I f
\left( y + \left( \lambda t_1 + (1 - \lambda) t_2 \right)  e_i, r \right) dr \\ \nonumber &
\leq \, R \int_I  \left[ \lambda f \left( y + t_1 e_i, r
\right) + (1-\lambda) f \left( y + t_2 e_i, r \right) \right] dr \\ & = R \left[ \lambda \pi(f)(y + t_1 e_i) + (1-\lambda)  \pi(f)(y + t_2 e_i) \right], \nonumber
\end{align}
and the lemma is proven. \qed

\medskip \textbf{Proof of Theorem \ref{thm_103}:}
We will prove by induction on the dimension $n$ that there exists a  transportation map $T$ from $f$ to $g$
such that
\begin{equation} \int_Q |T(x) - x|^2 f(x) dx \leq  \frac{40}{9} R^{2 } \left[ \int_Q S_T \{ g, f \} -
\left( \int_{Q} f \right) \log \frac{\int_Q g}{\int_Q f} \right].
\label{eq_1658} \end{equation}  The case $n=1$ is Proposition
\ref{prop_142}. Assume that the induction hypothesis was proven for dimension
$n-1$, and let us prove it for dimension $n$. Thus, suppose that we
are given a cube $Q \subset \RR^n$ and functions $f,g$ which satisfy
the assumptions of the theorem. In view of Lemma \ref{lem_2215}, we
may apply the induction hypothesis for
$$ \pi(Q), \pi(f), \pi(g). $$
Thus, there exists a transportation map $P: \overline{\pi(Q)}
\rightarrow \overline{\pi(Q)}$ from $\pi(f)$ to $\pi(g)$ such
that
\begin{align} \label{eq_1020}  \int_{\pi(Q)}  & |P(y) - y|^2 \pi(f)(y) dy \\ & \leq  \frac{40}{9} R^{2 }  \left[ \int_{\pi(Q)} S_P \{ \pi(g), \pi(f) \} -
\left( \int_{\pi(Q)} \pi(f) \right) \log \frac{\int_{\pi(Q)}
\pi(g)}{\int_{\pi(Q)} \pi(f)} \right]. \nonumber \end{align} For $y
\in \overline{\pi(Q)}$ let $T_y$ be the monotone transportation map
from $f_y$ to $g_{P(y)}$, a strictly-increasing  function which
is well-defined for almost any $y \in \pi(Q)$. We set
$$ T(y, r) = (P(y), T_y(r)) \quad \quad \quad \quad \text{for} \ (y,r) \in \overline{Q}. $$
Then $T$  transports the last coordinate monotonically. Hence,
according to Lemma \ref{lem_1232},
\begin{equation} \int_Q     \left |  T_y(r) -r \right |^2 f(y,r) dy dr   \, \leq \, \frac{40}{9} R^{2 } \left[ \int_Q S_{T} \{g,
f\} - \int_{\pi(Q)} S_{P} \{\pi(g), \pi(f)\} \right].
\label{eq_1021}
\end{equation}
It is straightforward to verify that the map $T$ is a transportation
map from $f$ to $g$. In fact,  the map $T$
is precisely  the Knothe
transportation map from \cite{knothe}. By summing (\ref{eq_1020})
and (\ref{eq_1021}), we conclude that
\begin{align}
\label{eq_1022} \int_Q  & \left[    \left|P(y) - y \right|^2 + \left
|  T_y(r) -r \right |^2 \right] f(y,r) dy dr  \\  \nonumber & \leq
\frac{40}{9} R^{2 } \left[ \int_Q S_{T} \{g, f \} - \left( \int_{\pi(Q)}
\pi(f) \right) \log \frac{\int_{\pi(Q)} \pi(g)}{\int_{\pi(Q)}
\pi(f)} \right]
\\ & =  \frac{40}{9} R^{2  } \left[ \int_Q S_{T} \{g, f\} - \left( \int_{Q} f \right) \log \frac{\int_{Q} g}{\int_{Q} f} \right]. \nonumber
\end{align}
All that remains is to note that when $x = (y,r)$,
$$ \left|T(x) - x \right|^2 =  \left | P(y) - y \right|^2 +     \left |  T_y(r) -r \right |^2. $$
The bound (\ref{eq_1658}) follows from (\ref{eq_1022}), and the
theorem is proven. \qed

\section{Log-concavity}
\label{sec4}

We begin this section with a  discussion of Definition \ref{tire}.
As we shall see, this definition fits very nicely with log-concave functions.
Given two probability densities $f$ and $g$ in  $\RR^n$ we write
$$ D( g \, || \, f ) = \int_{\RR^n}  \left[ \log \frac{g(y)}{f(y)} \right] g(y) dy $$
for the relative entropy or the Kullback–-Leibler divergence of $g$ from $f$.

\begin{lemma} Let $f, g: \RR^n \rightarrow [0, \infty)$ be probability densities.
Assume that $f$ is log-concave. Then,
$$  Tire( g \, || \, f ) \leq D( g \, || \, f ).
$$
\label{lem_1035}
\end{lemma}

\textbf{Proof:}  The function $f$ is differentiable almost-everywhere
 in the convex set $Supp(f)$ as it is a log-concave function. Denote $f =
e^{-\psi}$. From the convexity of $\psi$ we see that for any point
$x \in Supp(f)$ in which $f$ is differentiable,
$$ \psi(x) + \nabla \psi(x) \cdot (y-x) \leq \psi(y) \quad \quad \quad \quad (y \in Supp(f)). $$
Let $T$ be any transportation map from $f$ to $g$. Denoting $\vphi = -\log g$ we find that for almost any $x \in Supp(f)$,
\begin{align*}
 S_T \{g, f \}(x)   \, & = \,  \left[ f(x) \log \frac{g(T(x))}{f(x)} - \nabla f(x) \cdot
(T(x) - x) \right] \\ & = \, f(x)  \left[ \psi(x) - \vphi(T(x)) + \nabla \psi(x) \cdot (T(x)-x) \right]
\\ &  \leq \, f(x) \left[ \psi(T(x)) - \vphi(T(x)) \right] = f(x)
\log \frac{g(T(x))}{f(T(x))}.
\end{align*}
Since $T$ is a transportation map from $f$ to $g$, then by applying the change of variables $y = T(x)$
 we obtain
\begin{align*}
 \int_{Supp(f)} S_T \{ g, f \} &
\leq \int_{Supp(f)} \log \frac{g(T(x))}{f(T(x))} f(x) dx \\ & = \int_{Supp(g)} \left[ \log \frac{g(y)}{f(y)} \right]g(y) dy
= D(g \, || \, f) \end{align*}
and the lemma is proven. \qed

\medskip
For a log-concave density $f$, we may think about $Tire(g \, || \, f)$  as a parameter measuring the proximity of
$g$ to a translate of $f$.
Let us mention here additional upper bounds for $Tire(g \, || \, f)$. Let $f,g: \RR^n
\rightarrow [0, \infty)$ have finite, positive integrals and denote $\psi = -\log f$ and $\vphi = -\log
 g$. According to (\ref{eq_501}),
 \begin{equation} \sup_{y \in Supp(g)} S_{y} \{g, f \}(x) = \left[ \nabla f(x) \cdot x - f(x) \log f(x) \right] + f(x) \vphi^*( \nabla \psi(x)) \label{eq_1634} \end{equation}
where $\vphi^*(v) = \sup_{y \in Supp(g)} \left[ v \cdot y -
\vphi(y) \right]$ is the usual Legendre transform of  $\vphi$. Consequently, when $f$
is locally-Lipschitz and $x \mapsto |x| f(x)$ vanishes at infinity,
we have the bound
\begin{equation}  Tire(g \, || \, f) \leq \int_{\RR^n}
 \left[ \vphi^* \left( \nabla \psi(x) \right) \, - \, \log \left( \frac{\int g}{\int f} \cdot f(x) \right) \, - \, n \right] f(x) dx.
 \label{eq_956} \end{equation}
 Inequality (\ref{eq_956}) is perhaps less appealing than Lemma \ref{lem_1035}, yet it is applicable
 also in the non log-concave case.

\medskip Our original motivation for Definition \ref{tire}
is that at least in the smooth, log-concave case, the expression in (\ref{eq_1634}) equals
$
\left.
\partial h_{\eps}(x) \left/
\partial \eps \right. \right|_{\eps = 0}$ where
$$ h_{\eps}(x) = \sup_{y \in \RR^n} f(x + \eps y)^{1-\eps} g(x - (1
- \eps) y)^{\eps} \quad \quad \quad \quad (x \in \RR^n). $$
In other words, $Tire(g \, || \, f)$ is related to a kind of
 ``mixed volume'' of log-concave functions, see \cite[Section 3]{marg} for further explanations.

\medskip Suppose that
$\mu_1$ and $\mu_2$ are  Borel probability measures on $\RR^n$. The
{\it transportation cost} between $\mu_1$ and $\mu_2$ is defined to
be
$$ W_2^2(\mu_1, \mu_2) = \inf_{\gamma} \int_{\RR^n \times \RR^n}  |x - y|^2 d \gamma(x,y) $$
where the infimum runs over all couplings $\gamma$ of $\mu_1$ and
$\mu_2$, i.e., all Borel probability measures $\gamma$ on $\RR^n
\times \RR^n$ whose first marginal is $\mu_1$ and whose second
marginal is $\mu_2$. See, e.g., Villani's book \cite{villani} for
more information about the transportation metric $W_2$.
 The following theorem
reminds us of Talagrand's transportation-cost inequalities for
product measures from \cite{talagrand}.

\begin{theorem} Let $R \geq 1$ and let $Q \subset \RR^n$ be a unit cube  parallel to the axes. Suppose that $\mu$ is a probability measure on $Q$ with a log-concave density
$f$. Assume that
\begin{equation} f \left( \lambda x + (1 - \lambda) y \right) \leq R \left[ \lambda f(x) + (1-\lambda) f(y) \right]
\label{eq_1228} \end{equation} for any $0 < \lambda < 1$ and any $x,
y \in Q$ for which $x - y$ is proportional to one of the standard
unit vectors in $\RR^n$. Let $\nu$ be a probability measure on $Q$
that is absolutely continuous with respect to $\mu$. Then,
$$ W_2^2( \mu, \nu ) \leq C R^{2 } D \left( \nu \, || \, \mu \right) $$
where $D \left( \nu \, || \, \mu \right) = \int g (\log g) d \mu$
for $g = d \nu / d \mu$, the usual relative entropy functional, and
where $C \leq 40/9$ is a universal constant.
 \label{prop_1345}
\end{theorem}

\textbf{Proof:} By a standard approximation argument (e.g.,
convolve $\mu$ with a tiny gaussian and restrict to the cube $Q$), we
may assume that $f$ and $g$ are positive and $C^1$-smooth up to the boundary in $Q$, and in particular
both functions are positive and Lipschitz.
 According to Theorem \ref{thm_103} and Lemma  \ref{lem_1035},
 $$ W_2^2(\mu, \nu) \leq C R^{2 } \cdot Tire \left (g \, || \, f \right) \leq C R^2 \cdot D \left( \nu \, || \, \mu \right).
 $$

 \vspace{-24pt}  \qed

\medskip Transportation-cost inequalities such as Theorem \ref{prop_1345}
are the subject of the comprehensive survey by  Gozlan and L\'eonard
\cite{GZ}. The fact that transportation-cost inequalities  imply
concentration inequalities goes back to Marton \cite{marton}. The
following proof reproduces her argument, and is included here for
completeness.

\medskip \textbf{Proof of Theorem \ref{cor_1636}:}
Denote $E = Q \setminus (A + t B^n)$. If $\mu(E) = 0$ then there is
nothing to prove. Otherwise, we apply Theorem \ref{prop_1345} for the measure $\nu = \mu|_E$. Thus there exists
a coupling $\gamma$ of $\mu$ and $\mu|_E$  with
$$ \int_{Q \times E} \left |y - x \right|^2 d \gamma(x,y) \leq \frac{40}{9} R^{2} D( \nu \, || \, \mu) =  \frac{40}{9} R^{2} \cdot \log \frac{1}{\mu(E)}. $$
According to the Markov-Chebyshev inequality, there exists a subset $F \subseteq
Q \times E$ with $\gamma(F) \geq 41/81$ such that for any $(x,y) \in F$,
\begin{equation}
 \left| y - x \right|^2 \leq 9 R^{2  } \log \frac{1}{\mu(E)}.
\label{eq_1125}
\end{equation}
Since $\gamma$ is a coupling and $\mu(A) \geq 1/2$ with $\gamma(F) \geq 41/81$,
there exists $(x,y) \in F$ with $x \in A$. For such $(x,y)$,
$$ x \in A, y \in E \quad \quad \text{and} \quad \quad  \left| x - y \right| \leq 3  R \cdot \sqrt{\log \frac{1}{\mu(E)}} $$
where we used (\ref{eq_1125}).
However, all points in $E$ are of distance at least $t$ from all points of $A$. Consequently,
$$ t \leq 3 R \cdot \sqrt{\log \frac{1}{\mu(E)}}. $$
Therefore $\mu(E) \leq \exp(-t^2 / \alpha^2)$ for $\alpha = 3 R$
and $\mu(A + t B^n) \geq 1 - \exp(-t^2 / \alpha^2)$, as required. \qed

\medskip \textbf{Proof of Theorem \ref{thm_2140}:}
Let $T > 0$. Observe that the validity of both the assumptions and
the conclusions of the theorem is not altered under the scaling
$$ \ell \mapsto T \ell, \quad M \mapsto T^{-2} M. $$
We may thus normalize so that $\ell = 1$.
All that remains is to verify that the assumptions
of Theorem \ref{cor_1636} are satisfied with $R = e^{M/8}$.
Fix $i=1,\ldots,n$ and $x \in Q$ and denote $h(t) = \psi(x + t e_i)$.
Then $h$ is well-defined on a certain interval $I \subset \RR$ of length one, and our goal is to show
that for any $a, b \in I$ and $0 < \lambda < 1$,
\begin{equation}  e^{- h(\lambda a + (1 - \lambda) b)} \leq e^{M/8} \left[ \lambda e^{-h(a)} + (1 -
\lambda) e^{-h(b)} \right]. \label{eq_2127}
\end{equation}
In view of the arithmetic/geometric means inequality, the desired
inequality (\ref{eq_2127}) would follow once we establish that
\begin{equation}   - h(\lambda a + (1 - \lambda) b) \leq M / 8
- \lambda h(a) - (1 -
\lambda) h(b). \label{eq_1731}
\end{equation}
In order to prove (\ref{eq_1731}), we use our assumption that $ h^{\prime \prime}(t) \leq M$
in the interval $I$.
According to the Taylor theorem, for any $x, y \in I$,
\begin{equation}  \label{eq_239} h(y) - h(x) - h^{\prime}(x) (y - x) \leq M \frac{(x-y)^2}{2}.
\end{equation} We will apply
inequality  (\ref{eq_239}) for $y=a,b$ and $x = \lambda a + (1 -
\lambda) b$, then add the resulting inequalities with coefficients
$\lambda$ and $1 -\lambda$. This yields
\begin{equation}  \lambda h(a) + (1 -
\lambda) h(b) - h(\lambda a + (1 - \lambda) b) \leq M \frac{
\lambda  (1 - \lambda) (b-a)^2  }{2} \leq \frac{M}{8} \label{eq_159}
\end{equation}
as $\lambda (1 - \lambda) \leq 1/4$ and $|b-a| \leq 1$. The inequality (\ref{eq_1731})  follows from (\ref{eq_159}). \qed

\medskip It is well-known (see, e.g., V. Milman and Schechtman \cite[Section 2 and Appendix V]{MS}) that Theorem \ref{thm_2140}
implies a concentration inequality for Lipschitz functions as follows:

\begin{corollary} Let $Q, \mu, \alpha$ be as in Theorem \ref{thm_2140}
(or as in Theorem \ref{cor_1636}).
Let $f: Q \rightarrow \RR$ be a $1$-Lipschitz function, i.e., $|f(x) - f(y)| \leq |x-y|$
for any $x,y \in Q$. Denote $E = \int_Q f d\mu$. Then, for any $t > 0$,
$$ \mu \left \{ x \in Q \, ; \, \left| f(x) - E \right| \geq t \right \} \leq C e^{-c t^2 / \alpha^2} $$
where $c, C > 0$ are universal constants.
\label{cor_1329}
\end{corollary}

In particular, we deduce from Corollary \ref{cor_1329} that in the notation of Theorem \ref{thm_2140},
\begin{equation}  Cov(\mu) \leq C \alpha^2 \cdot Id \label{eq_1330} \end{equation}
in the sense of symmetric matrices, where $Cov(\mu)$ is the covariance matrix of the probability measure $\mu$
and $C > 0$ is a universal constant.

\begin{remark} {\rm Regarding the dependence of $\alpha(\ell,M)$ on $M$ in Theorem \ref{thm_2140}:
Let $X_0,\ldots,X_n$ be independent standard Gaussian random variables. Consider the random vector
$$ Y = \frac{(X_1,\ldots,X_n)}{100 \sqrt{\log n}} +  \frac{(X_0,\ldots,X_0)}{100 \sqrt{\log n}}, $$
and let $Z$ be the conditioning of $Y$ to the cube $Q =
[-1/2,1/2]^n$. Denote by $\mu$ the distribution of $Z$, a
probability measure on $Q$. It is not too  difficult to verify that
$\mu$ satisfies the requirements of Theorem \ref{thm_2140}
 with $\ell=1$ and $M = C \log n$. Set $A = \{ x \in Q ; \sum_i x_i \leq 0 \}$. Then $\mu(A) = 1/2$.
 However, one may compute that for any $t \leq c n^{1/2} / \sqrt{\log n}$,
 $$ \mu(A + t B^n) \leq 2/3. $$
This shows that $ \alpha(1, C \log n) \geq c n^{1/2} / \sqrt{\log
n}$. Therefore the exponential dependence of the dimension-free
expression $\alpha(\ell,M)$ on $\ell^2 M$ is inevitable. A simple
variant of this example shows that it is also impossible to replace
the cube $Q$ of sidelength $\ell$ in Theorem \ref{thm_2140} by a
Euclidean ball of radius $\ell \sqrt{n}$. For another example in
which the cube behaves better than the Euclidean ball, see
\cite[Corollary 3]{ptrf}. } \label{rem_1703}
\end{remark}

\medskip It was explained by E. Milman \cite{emanuel}
that in the log-concave case, Gaussian concentration inequalities,
quadratic transportation-cost inequalities, and log-Sobolev
inequalities are all essentially equivalent up to universal
constants. In particular, by using the results of Otto and Villani
\cite[Corollary 3.1]{OV}, we deduce from Theorem \ref{prop_1345} the
following log-Sobolev and Poincar\'e inequalities:

\begin{corollary} Let $\ell, M, Q, \mu$ be as in Theorem \ref{thm_2140}
(or as in Theorem \ref{cor_1636}, with $\ell = 1$ and $R = e^{M/8}$).
 Then,  for any locally-Lipschitz function $f: Q \rightarrow \RR$ with $\int_Q f^2 d\mu = 1$,
 $$ \int_Q f^2 \log \left( f^2 \right) d \mu \leq C_1 \ell^2 e^{M \ell^2/4} \int_Q |\nabla f|^2 d \mu, $$
 and for any integrable, locally-Lipschitz function $f: Q \rightarrow \RR$ with $\int_Q f d\mu = 0$,
 $$ \int_Q f^2 d \mu \leq C_2 \ell^2 e^{M \ell^2/4} \int_Q |\nabla f|^2 d \mu. $$
 Here, $C_1 \leq 160/9$ and $C_2 \leq 20/9$ are universal constants. \label{cor_1509}
 \end{corollary}

It is conceivable that Theorem \ref{thm_2140} and Corollary
\ref{cor_1509} will turn out to be relevant to the analysis of
lattice models in physics. For instance, one may suggest an Ising
model with bounded, real spins as in Royer \cite[Section 4.2]{royer}
in which the assumptions of Theorem \ref{thm_2140} are satisfied.
Essentially, we require that the spins lie in the interval $[-1,1]$,
that the entire Hamiltonian is convex (just convex, not
strictly-convex) and that the second derivatives of the pairwise
potentials and the self-interactions are bounded. Perhaps the
logarithmic Sobolev inequality of Corollary \ref{cor_1509} may be
of some use in this context.

\section{Yet another approach for Theorem \ref{prop_1345}}

In this section we present a sketch of an alternative proof of Theorem \ref{prop_1345},
in spirit of the transportation arguments of Cordero-Erausquin
\cite{DCE}. The derivation below is applicable for the two types of transportation
maps, Brenier and Knothe.

\medskip Let $f, g, \mu$ and $\nu$ satisfy the assumptions of Theorem \ref{prop_1345}.
As is explained above, it suffices to consider the case where $f$ and $g$
are positive, Lipschitz functions. In particular, it is well-known
that both the Brenier map and the Knothe map from $\mu$ to $\nu$ are $C^1$-smooth
up to the boundary (see Cordero-Erausquin \cite{C2}).

\medskip
Denote $\psi = -\log f$, a convex function.
Let $F$ be any smooth transportation map from $\mu$ to $\nu$. Then, similarly to (\ref{eq_1047}) above,
we have
$$
 \log \left| \det F^{\prime}(x) \right| = -\psi(x) + \psi(F(x)) - \log g(F(x)) \quad \quad \quad (x \in Q)
 $$ where $F^{\prime}(x)$ is the $n \times n$ matrix
which is the derivative of $F$. In the case where $F$ is the Brenier map, the matrix $F^{\prime}(x)$
is symmetric and positive-definite. In the case where $F$ is the Knothe map, the matrix $F^{\prime}(x)$ is upper-triangular with positive entries on the diagonal.
In both cases, denoting $F = (F_1,\ldots,F_n)$,
\begin{equation}
 \log \left| \det F^{\prime}(x) \right| = \log \det F^{\prime}(x) \leq \sum_{i=1}^n \log \partial^i F_i(x) \quad \quad \quad (x \in Q). \label{eq_1437}
\end{equation}
Indeed, in the Knothe case (\ref{eq_1437}) is simply an equality, while in the Brenier case we may use Hadamard's determinant inequality
in order to establish (\ref{eq_1437}).
Next, denote $\theta(x) = F(x) - x$, so that $\partial^i F_i(x) = 1 + \partial^i \theta_i(x)$.
We use the elementary inequality for the logarithm function in (\ref{eq_1124}) and obtain
$$
\sum_{i=1}^n \log
\partial^i F_i(x) \leq \sum_{i=1}^n \left[ \partial^i \theta_i(x)  -\frac{3}{10} \Lambda \left( \partial^i \theta_i(x) \right) \right].
$$
The convexity of $\psi$ implies that
$$ \psi(F(x)) - \psi(x) \geq \nabla \psi(x) \cdot (F(x) - x) = \nabla \psi(x) \cdot \theta(x) = \sum_{i=1}^n \theta_i(x) \partial^i \psi(x) . $$
Combining all of the above, we arrive at the inequality
\begin{equation}
\log g(F(x))
\, \geq \, \frac{3}{10} \sum_{i=1}^n   \Lambda \left( \partial^i \theta_i(x) \right)
\, - \,
\sum_{i=1}^n \left[ \partial^i \theta_i(x) -  \partial^i \psi(x) \cdot \theta_i(x) \right], \label{eq_1508}
\end{equation}
valid pointwise in $Q$. Here comes a fundamental property of both the Brenier map and the Knothe map:
In both cases, the map $F$ preserves each of the $(n-1)$-dimensional facets of the cube $Q$. In other words,
let $A_1^{\pm}, \ldots, A_n^{\pm}$ be an enumeration of all of the $2n$ facets of dimension $n-1$
of the cube $Q$. Assume that $\pm e_i$ is the outer unit normal to the cube $Q$
at the facet $A_i^{\pm}$. We claim that  for any $i$,
\begin{equation}
x \in A_i^{\pm} \quad \quad  \Longrightarrow \quad \quad \theta_i(x) = 0. \label{eq_1459} \end{equation}
It is quite clear that (\ref{eq_1459}) holds in the case of the Knothe map.
In order to argue for (\ref{eq_1459}) in the Brenier case, recall that here
$$ (F(x) - F(y)) \cdot (x -y) > 0 \quad \quad \quad \quad (x,y \in Q, \ x \neq y) $$
as $F$ is the gradient of a strictly-convex function. In particular, when $F(x) \in A_i^{\pm}$ then necessarily
$x \pm t e_i \not \in Q$ for $t \in (0, \eps)$ for some $\eps > 0$. Hence $F(x) \in A_i^{\pm}$ implies that $x \in A_i^{\pm}$.
Arguing similarly for the inverse map $F^{-1}$, which is the Brenier map from $\nu$ to $\mu$, we conclude that
(\ref{eq_1459}) holds true in the Brenier case as well.

\medskip
We may now multiply (\ref{eq_1508}) by $e^{-\psi}$ and integrate over the cube $Q$. Observe that for any $i=1,\ldots,n$,
$$ \int_{Q} \left[ \partial^i \theta_i(x) -  \partial^i \psi(x) \cdot \theta_i(x) \right] e^{-\psi} =
\int_{Q} \partial^i \left(  \theta_i e^{-\psi} \right)  = 0 $$
thanks to the boundary condition (\ref{eq_1459}). Furthermore, this boundary condition 
allows us to use  the one-dimensional  Lemma \ref{lem_233}, and conclude that
$$  \int_{Q} \Lambda \left(\theta_i(x) \right) d \mu(x) \leq \frac{4 R^2}{3} \int_{Q} \Lambda \left( \partial^i \theta_i(x) \right) d \mu(x)  $$
for $i=1,\ldots,n$. We therefore obtain
$$
\int_{Q} \sum_{i=1}^n \Lambda \left(\theta_i(x) \right) d \mu(x) \leq \frac{40}{9} R^2 \int_{Q} \left[ \log g(F(x)) \right] d \mu(x)
= \frac{40}{9} R^2 \int_{Q} \left[ \log g(y) \right] d \nu(y).
$$
It remains to note that always $|\theta_i(x)| = |F_i(x) - x_i| \leq 1$ since $Q$ is a unit cube. Consequently $\Lambda(\theta_i(x)) = |\theta_i(x)|^2$ and hence,
$$ \int_{Q} |F(x) - x|^2 d \mu(x) = \int_{Q} |\theta(x)|^2 d \mu(x) \leq \frac{40}{9} R^2 \int_{Q} \left[ \log g(y) \right] d \nu(y) = \frac{40}{9} R^2 \cdot D \left( \nu \, || \, \mu \right). $$
This finishes the sketch of the alternative proof of Theorem \ref{prop_1345}.\qed

{
}


\begin{thebibliography}{99}
\bibitem{BE} Bakry, D., \'Emery, M., {\it Diffusions hypercontractives.}  S\'eminaire de probabilit\'es, XIX, 1983/84,
Lecture Notes in Math., Vol. 1123, Springer, Berlin, (1985),
177--206.

\bibitem{barthe} Barthe, F., {\it In\'egalit\'es de Brascamp-Lieb et convexit\'e}. C. R. Acad. Sci. Paris
S\'er. I Math., Vol. 324., No. 8, (1997),  885--888.

\bibitem{barthe2} Barthe, F., {\it On a reverse form of the Brascamp-Lieb inequality. }
Invent. Math., Vol. 134, No. 2, (1998),  335--361.

\bibitem{brenier}  Brenier, Y., {\it Polar factorization and monotone rearrangement of vector-valued functions. }
Comm. Pure Appl. Math., Vol. 44, No. 4, (1991), 375--417.


\bibitem{cheeger}  Cheeger, J., {\it A lower bound for the smallest eigenvalue of the Laplacian. }
Problems in analysis (Papers dedicated to Salomon Bochner, 1969),
Princeton Univ. Press, Princeton, NJ, (1970), 195--199.

\bibitem{DCE} Cordero-Erausquin, D., {\it Some applications of mass transport to Gaussian-type inequalities.}
Arch. Ration. Mech. Anal., Vol. 161, No. 3, (2002), 257–-269.

\bibitem{C2} Cordero-Erausquin, D., {\it Sur le transport de mesures p\'eriodiques.}
 C. R. Acad. Sci. Paris S\'er. I Math., Vol. 329, No. 3, (1999),  199–-202.

\bibitem{DZ} Dembo, A., Zeitouni, O., {\it
 Transportation approach to some concentration inequalities in product spaces.}
  Electron. Comm. Probab., Vol. 1, No. 9, (1996), 83--90.

\bibitem{EK} Eldan, R., Klartag, B., {\it Dimensionality and the stability of the Brunn-Minkowski inequality.}
To appear in Ann. Sc. Norm. Super. Pisa. Available under \verb"http://arxiv.org/abs/1110.6584"

\bibitem{GZ} Gozlan, N., L\'eonard, C., {\it Transport inequalities. A survey. }
Markov Processes and Related Fields, Vol. 16, (2010), 635--736.

\bibitem{HM} Henstock, R., Macbeath, A. M., {\it On the measure of sum-sets. I. The
theorems of Brunn, Minkowski, and Lusternik. } Proc. London Math.
Soc. (3), Vol. 3, (1953), 182--194.


\bibitem{marg} Klartag, B., {\it
Marginals of geometric inequalities}.
Geometric Aspects of Functional Analysis, Lecture Notes in Math., Vol. 1910, Springer (2007), 133--166.

\bibitem{ptrf} Klartag, B., {\it A Berry-Esseen type inequality for convex bodies with an
unconditional basis. } Probab. Theory Related Fields, Vol. 145, No.
1-2, (2009), 1--33.

\bibitem{knothe} Knothe, H., {\it Contributions to the theory of convex
bodies.} Michigan Math. J., Vol. 4, (1957), 39--52.

\bibitem{ledoux_product} Ledoux, M., {\it On Talagrand's deviation inequalities for product measures}.
ESAIM: Probab. Statist., Vol. 1, (1997), 63--87.

\bibitem{leoni} Leoni, G., {\it A first course in Sobolev spaces.}
Graduate Studies in Mathematics, 105. Amer. Math. Soc., Providence, RI, 2009.

\bibitem{marton} Marton, K., {\it A measure concentration inequality for contracting Markov
chains.} Geom. Funct. Anal., Vol. 6, No. 3, (1996), 556--571.



\bibitem{mccann} McCann, R. J., {\it A convexity principle for interacting gases. }
Adv. Math., Vol. 128, No. 1, (1997), 153--179.


\bibitem{emanuel} Milman, E., {\it
Isoperimetric and concentration inequalities - equivalence under
curvature lower bound}. Duke Math. J., Vol. 154, No. 2, (2010),
207--239.

\bibitem{MS} Milman, V. D., Schechtman, G., {\it Asymptotic theory of finite-dimensional normed spaces. }
Lecture Notes in Math., Vol. 1200, Springer, Berlin, 1986.

\bibitem{OV} Otto, F.,  Villani, C., {\it Generalization of an inequality by Talagrand and links with the logarithmic
Sobolev inequality. } J. Funct. Anal., Vol. 173, No. 2, (2000),
361--400.

\bibitem{royer} Royer, G., {\it An initiation to logarithmic Sobolev inequalities.}
 SMF/AMS texts and monographs, 14. American Mathematical Society, Providence, RI;
Soci\'et\'e Math\'ematique de France, Paris, 2007.

\bibitem{talagrand} Talagrand, M., {\it Transportation cost for Gaussian and other
product measures. } Geom. Funct. Anal., Vol. 6, No. 3, (1996),
587--600.

\bibitem{talagrand_ihes} Talagrand, M., {\it  Concentration of measure and isoperimetric inequalities in product spaces.}
 Inst. Hautes \'Etudes Sci. Publ. Math., No. 81, (1995), 73--205.

\bibitem{villani} Villani, C., {\it Topics in optimal transportation. }
Graduate Studies in Mathematics, 58. Amer. Math. Soc.,
Providence, RI, 2003.

\end{thebibliography}
\end{document}